\documentclass[aos,seceqn,citesort,dvips]{arximspdf}
\usepackage{graphics}
\doi{10.1214/09-AOS755}
\volume{38}
\issue{3}
\pubyear{2010}
\firstpage{1433}
\lastpage{1459}

\makeatletter

\newproclaim{example}{Example}

\newtheorem{theorem}{Theorem}

\newproclaim{definition}{Definition}

\makeatother

\begin{document}
\begin{frontmatter}

\title{Optional P\'{o}lya tree and Bayesian inference}
\runtitle{Optional P\'{o}lya tree}

\pdftitle{Optional Polya tree and Bayesian inference}

\begin{aug}
\author[A]{\fnms{Wing H.} \snm{Wong}\corref{}\thanksref{t1}\ead[label=e1]{whwong@stanford.edu}} and
\author[A]{\fnms{Li} \snm{Ma}\thanksref{t2}\ead[label=e2]{ma2@stanford.edu}}
\runauthor{W. H. Wong and L. Ma}
\affiliation{Stanford University}
\address[A]{Department of Statistics\\
Stanford University\\
390 Serra Mall\\
Stanford, California 94305\\
USA\\
\printead{e1}\\
\phantom{E-mail: }\printead*{e2}} 
\end{aug}

\thankstext{t1}{Supported in part by NSF Grants DMS-05-05732 and DMS-09-06044.}
\thankstext{t2}{Supported by a Gerhard Casper Stanford Graduate Fellowship.}

\received{\smonth{5} \syear{2009}}
\revised{\smonth{9} \syear{2009}}

%
\begin{abstract}
We introduce an extension of the P\'{o}lya tree approach for
constructing distributions on the space of probability
measures. By using optional stopping and optional choice of
splitting variables, the construction gives rise to
random measures that are absolutely continuous with piecewise smooth
densities on partitions that can adapt to fit the data. The
resulting ``optional P\'{o}lya tree'' distribution has
large support in total variation topology and yields posterior
distributions that are also optional P\'{o}lya trees with computable
parameter values.
\end{abstract}

%

\pdfkeywords{62F15, 62G99, 62G07, Polya tree,
Bayesian inference, nonparametric, recursive partition,
density estimation}

\end{frontmatter}

\section{Introduction}\label{sec1}

Ferguson \cite{ferg73} formulated two criteria for desirable prior
distributions on the space of probability measures: (i) The support of
the prior should be large with respect to a suitable topology, and (ii)
the corresponding posterior distribution should be analytically
manageable. Extending the work by Freedman \cite{free63} and Fabius
\cite{fabius64},
he introduced the Dirichlet process as a prior that
satisfies these criteria. Specifically, assuming for simplicity that
the parameter space $\Omega$ is a bounded interval of real numbers, and
the base measure in the Dirichlet process prior is the Lebesgue
measure, then the prior will have positive probability in all weak
neighborhoods of any absolutely continuous probability measure, and
given i.i.d. observations, the posterior distribution is also a
Dirichlet process with its base measure obtainable from that of the
prior by the addition of delta masses at the observed data points.

While these properties made it an attractive prior in many Bayesian
nonparametric problems, the use of the Dirichlet process prior is
limited by its inability to generate absolutely continuous
distributions; that is, a random probability measure sampled from the
Dirichlet process prior is almost surely a discrete measure
\mbox{\cite{black73,blackmac73,ferg73}}. Thus in applications that require the
existence of
densities under the prior, such as the estimation of a density from a
sample \cite{lo84} or the modeling of error distributions in location
or regression problems \cite{dia86}, there is a need for alternative
ways to specify the prior. Lo \cite{lo84} proposed an elegant prior in
the space of densities by assuming the density is a mixture of kernel
functions where the mixing distribution is modeled by a Dirichlet
process. Under Lo's model, the random distributions are guaranteed to
have smooth densities and the predictive density is still analytically
tractable. However the degree of smoothness is not adaptive.

Another approach to deal with the discreteness problem is
to use P\'{o}lya tree priors \cite{ferg74}. This class of random
probability measures includes the Dirichlet process as a special case
and yet is itself a special case of the more general class of ``tail
free'' processes previously studied by Freedman \cite{free63}. P\'{o}lya
tree prior
satisfies Ferguson's two criteria. First, it is possible to construct
P\'{o}lya tree priors with positive probability in neighborhoods around
arbitrary positive densities \cite{lav92}. Second, the posterior
distribution arising from a P\'{o}lya tree prior is available in close form
\cite{ferg74}. Further properties and applications of P\'{o}lya tree
priors are found in \cite{hutter09,lav94,ghosh03,hans06} and \cite{maul92}.

In this paper we study the extension of the P\'{o}lya tree prior
construction by allowing
optional stopping and randomized partitioning schemes. To motivate
optional stopping, consider the standard construction of the P\'{o}lya tree
prior for probability measures in an interval $\Omega$. The interval
is recursively bisected into subintervals. At each stage, the
probability mass already assigned to an interval is randomly divided
and assigned into its subintervals according to the independent draw
of a Beta variable. However, in order for the prior to generate
absolutely continuous measures, it is necessary for the parameters in
the Beta distribution to increase rapidly as the depth of the
bisection increases, that is, as we move into more and more refined
levels of partitioning \cite{kraft64}.

In any case, even when the construction yields a random distribution
with density, with probability $1$ the density will have discontinuity
almost everywhere. The use of Beta variables with large magnitudes
for its parameters, although useful in forcing the random distribution
to be absolutely continuous, has the effect of severely constraining
our ability to allocate conditional probability to represent
faithfully the data distributions within small intervals. To resolve
this conflict between smoothness and faithfulness to the data
distribution, one can introduce an optional stopping variable for each
subregion obtained in the partitioning process~\cite{hutter09}. By
putting uniform
distributions within each stopped subregion, we can achieve the goal
of generating absolutely continuous distributions without having to
force the Beta parameters to increase rapidly. In fact, we will be
able to use Jeffrey's rule of Beta ($\frac12,\frac12$) in the
inference of conditional probabilities, regardless of the depth of the
subregion in the partition tree. We believe this is a desirable
consequence of optional stopping.

Our second extension is to allow randomized partitioning. Standard
P\'{o}lya tree construction relies on a fixed scheme for partitioning. For
example in \cite{hans06} a $k$-dimen\-sional rectangle is recursively
partitioned where in each stage of the recursion the subregions are
further divided into $2^k$ quadrants by bisecting each of the $k$
coordinate variables. In contrast, when recursive partitioning is used
in other statistical problems, it is customary to allow flexible
choices of the variables to use to further divide a subregion. This
allows the subregion to take very different shapes depending on the
information in the data. The data-adaptive nature of the recursive
partitioning is a reason for the success of tree-based learning
methodologies such as CART \cite{brei84}. Thus it is desirable to
allow P\'{o}lya tree priors to use partitions that are the result of
randomized choices of divisions in each of the subregions at each
stage of the recursion. Once the partitioning is randomized in the
prior, the posterior distribution will give more weights on those
partitions that provide better fits to the data. In this way the data
is allowed to influence the choice of the partitioning. This will be
especially useful in high-dimensional applications.

In Section \ref{sec2} we introduce the construction of ``Optional P\'{o}lya trees''
that allow optional stopping and randomized partitioning. It is shown
that this construction leads to priors that give absolutely continuous
distributions almost surely. We also show how to specify the prior so
that it has positive probability in all total variation neighborhoods
in the space of absolutely continuous distributions on $\Omega$. In
Section~\ref{sec3} we show that the use of optional P\'{o}lya tree priors will
lead to posterior distributions that are also optional P\'{o}lya trees. We
present a recursive algorithm for the computation of the parameters
governing the posterior optional P\'{o}lya tree. These results ensure that
Ferguson's two criteria are satisfied by optional P\'{o}lya tree priors,
but now on the space of absolutely continuous probability measures.
In this section, we also show that the posterior P\'{o}lya tree
is weakly consistent in the sense that asymptotically it concentrates
all its probability in any weak neighborhood of a true distribution
whose density is bounded. In Section \ref{sec:den}, we develop and test the
optional P\'{o}lya
tree approach to density estimation in Euclidean space. Concluding
remarks are given in Section \ref{sec:concl}.

We end this introduction with brief remarks on related works. The
important idea of early stopping was first introduced by
Hutter \cite{hutter09}. Ways to attenuate the dependency of P\'{o}lya trees
on the
partition include mixing the base measure used to define the tree
\cite{lav92,lav94,hans06}, random perturbation of the
dividing boundary in the partition of intervals \cite{paddock03} and
the use of positively correlated variables for the
conditional probabilities at each level of the tree definition
(Nieto-Barajas and M\"{u}ller~\cite{BarajasandMuller2009}).
Compared to these works, our approach allows not only
early stopping but also randomized choices of the splitting
variables. This provides a much richer class of partitions than
previous models and raises the new challenge of learning the partition
based on the observed data. We show that under mild conditions
such learning is achievable by finite computation. We also
provide a relatively complete mathematical foundation which
represents the first theory for Bayesian density estimation based on
recursive partitioning. Although a Bayesian version of recursive
partitioning has been proposed previously (Bayesian CART,
\cite{denison98}), it was formulated for a different problem (classification
instead of density estimation). Furthermore, it studied mainly model
specification and computational algorithm, and did not discuss the
mathematical and asymptotic properties of the method.\looseness=1

\section{Optional P\'{o}lya tree}\label{sec2}

We are interested in constructing random probability measures on a
space $(\Omega,\mu)$. $\Omega$ is either finite or a bounded rectangle
in $\mathbb{R}^p$. In this paper we assume for simplicity that $\mu$
is the
counting measure in the finite case and the Lebesgue measure in the
continuous case. Suppose that $\Omega$ can be partitioned in $M$
different ways; that is, for $j=1,2,\ldots,M$,
\[
\Omega=\bigcup_{k=1}^{K^j}\Omega_k^j \qquad\mbox{where
}\Omega_k^j\mbox{'s are disjoint}.
\]
Each $\Omega_k^j$, called a level-1 elementary region, can in turn be
divided into level-2 elementary regions. Assume there are
$M_{k_1}^{j_1}$ ways to divide $\Omega_{k_1}^{j_1}$; then for
$j_2=1,\ldots,M_{k_1}^{j_1}$, we have
\[
\Omega_{k_1}^{j_1}=\bigcup_{k_2=1}^{K_{k_1}^{j_1j_2}}\Omega
_{k_1k_2}^{j_1j_2}.
\]
In general, for any level-$k$ elementary region $A$, we assume there
are $M(A)$ ways to partition it; that is, for $j=1,2,\ldots,M(A)$,
\[
A=\bigcup_{k=1}^{K^j(A)}A_k^j.
\]
Let $\mathcal{A}^k$ be the set of all possible level-$k$ elementary regions,
and $\mathcal{A}^{(k)}=\bigcup_{l=1}^k\mathcal{A}^l$. If $\Omega$
is finite, we assume
that $\mathcal{A}^k$ separates points in $\Omega$ if $k$ is large
enough. If
$\Omega$ is a rectangle in $\mathbb{R}^p$, we assume that every open set
$B\subset\Omega$ is approximated by unions of sets in $\mathcal
{A}^{(n)}$, that is,
$\exists B_n\uparrow B$ where $B_n$ is a finite union of
disjoint regions in~$\mathcal{A}^{(n)}$.
\begin{example}\label{ex1}
\begin{eqnarray*}
\Omega &=& \bigl\{x=(x_1,\ldots,x_p)\dvtx x_i\in\{1,2\} \bigr\},\\
\Omega_k^j &=& \{x\dvtx x_j=k\},\qquad k=1\mbox{ or }2,\\
\Omega_{k_1k_2}^{j_1j_2} &=& \{x\dvtx x_{j_1}=k_1, x_{j_2}=k_2\},\qquad\mbox{ etc.}
\end{eqnarray*}
In this example, the number of ways to partition a level-$k$
elementary region decreases as $k$ increases.
\end{example}
\begin{example}\label{ex2}
\[
\Omega= \{(x_1,x_2,\ldots,x_p)\dvtx x_i\in[0,1] \}\subset\mathbb{R}^p.
\]
If $A$ is a level-$k$ elementary region (a rectangle), and $m_j(A)$ is
the midpoint of the range of $x_j$ for $A$, we set $A_1^j=\{x \in
A\dvtx x_j\leq m_j(A)\}$ and $A_2^j=A\setminus A_1^j$. There are exactly
$M(A)=p$ ways
to partition each $A$, regardless of its level.
\end{example}

Once a system to generate partitions has been specified as above, we
can formally define recursive partitions as follows. A recursive
partition of depth $k$ is a series of decisions
$J^{(k)}=(J_1,J_2,\ldots,J_k)$ where $J_l$ represents all the
decisions made at level $l$ to decide, for each region produced at the
previous level, whether or not to stop partitioning it further and if
not, which way to use to partition it. Once we have decided not
to partition a region, then it will remain intact at all subsequent
levels. Thus each $J^{(k)}$ specifies a partition of $\Omega$ into a
subset of regions in~$\mathcal{A}^{(k)}$.

We use a recursive procedure to produce a random recursive partition
of $\Omega$ and a random probability measure $Q$ that is uniformly
distributed within each part of the partition. Suppose after $k$ steps
of the recursion, we have obtained a random recursive partition
$\mathbf{J}^{(k)}$ and we write
\[
\Omega=T_0^k\cup T_1^k,
\]
where
\begin{eqnarray*}
T_0^k&=&\bigcup_{i=1}^IA_i \qquad\mbox{is a union of disjoint }A_i\in
\mathcal{A}
^{(k-1)},\\
T_1^k&=&\bigcup_{i=1}^{I'}A'_i \qquad\mbox{is a union of disjoint }A'_i\in
\mathcal{A}^k.
\end{eqnarray*}
The set $T_0$ represents the part of $\Omega$ where the partitioning
has already been stopped and $T_1$ represents the complement. In
addition, we have also obtained a random probability measure $Q^{(k)}$
on $\Omega$ which is uniformly distributed within each region in
$T_0^k$ and $T_1^k$.

In the $(k+1)$th step, we define $Q^{(k+1)}$ by further partitioning
of the regions in $T_1^k$ as follows. For each elementary region $A$
in the above decomposition of $T_1^k$, generate an independent random variable,
\[
S\sim\operatorname{Bernoulli} (\rho).
\]
If $S=1$, stop further partitioning of $A$ and add it to the set of
stopped regions. If $S=0$, draw $J\in\{1,2,\ldots,M(A)\}$ according to
a nonrandom vector $\bolds\lambda(A)=(\lambda_1,\ldots,\lambda_{M(A)})$,
called the selection probability vector, that is,
$P(J=j)=\lambda_j$ and $\sum_{l=1}^{M(A)}\lambda_l=1$.
If $J=j$, apply the $j$th way of partitioning $A$,
\[
A=\bigcup_{l=1}^KA_l^j \qquad\mbox{(here $K$ depends on $A$ and $j$)}
\]
and set $Q^{(k+1)}(A_l^j)=Q^{(k)}(A)\theta_l^j$ where
$\bolds\theta^j=(\theta_1^j,\ldots,\theta_K^j)$ is generated from a
Dirichlet distribution with parameter $(\alpha_1^j,\ldots,\alpha_K^j)$.
The nonrandom vector $\bolds\alpha^j=\bolds\alpha^j(A)$ is
referred to as
the assignment weight vector.

After this step, we have obtained $T_0^{k+1}$ and $T_1^{k+1}$, the
respective unions of the stopped and continuing regions. Clearly
\begin{eqnarray*}
\Omega &=& T_0^{k+1}\cup T_1^{k+1},\\
T_0^{k+1} &\supset& T_0^k,\qquad T_1^{k+1}\subset T_1^k.
\end{eqnarray*}
The new measure $Q^{(k+1)}$ is then defined as a refinement of
$Q^{(k)}$. For $B\subset T_0^{(k+1)}$, we set
\[
Q^{(k+1)}(B)=Q^{(k)}(B).
\]
For $B\subset T_1^{(k+1)}$ where $T_1^{k+1}$ is partitioned as
\[
T_1^{k+1}=\bigcup_{i=1}^JA_i,\qquad A_i\in\mathcal{A}^{k+1},
\]
we set
\[
Q^{(k+1)}(B)=\sum_{i=1}^JQ^{(k+1)}(A_i) \biggl(\frac{\mu(A_i\cap B)}{\mu
(A_i)} \biggr).
\]
Recall that for each $A_i$ in the partition of $T_1^{k+1}$, we have
already generated its $Q^{(k+1)}$ probability.

Let $\mathcal{F}^{(k)}$ be the $\sigma$-field of events generated by
all random
variables used in the first $k$ steps; the stopping probability
$\rho=\rho(A)$ is required to be measurable with respect to
$\mathcal{F}^{(k)}$. The specification of $\rho(\cdot)$ is called
the stopping
rule. In this paper we are mostly interested in the case when
$\rho(\cdot)$ is an ``independent stopping rule;'' that is, $\rho
(A)$ is
a pre-specified constant for each possible elementary region $A$.
However in some applications it is useful to let $\rho(A)$ depend on
$Q^{(k)}(A)$.

Let $\mathcal{A}^{(\infty)}=\bigcup_{k=1}^\infty\mathcal{A}^k$ be
the set of all
possible elementary regions.
\begin{theorem}\label{thm1}
Suppose there is a $\delta>0$ such that with probability $1$,
$1-\delta>\rho(A)>\delta$ for any region $A$ generated during any step in the
recursive partitioning process. Then with probability $1$, $Q^{(k)}$
converges in variational distance to a probability measure $Q$ that is
absolutely continuous with respect to $\mu$.
\end{theorem}
\begin{definition}
The random probability measure $Q$ defined in Theorem \ref{thm1} is
said to have an optional P\'{o}lya tree distribution with parameters
$\bolds\lambda, \bolds\alpha$ and stopping rule $\rho$.
\end{definition}
\begin{pf*}{Proof of Theorem \protect\ref{thm1}}
We only need to prove this for the case when $\Omega$ is a bounded
rectangle. We can think of $Q^{(k)}$'s as being generated in two
steps.
\begin{enumerate}
\item Generate the nonstopped version $Q^{*(k)}$ by recursively
choosing the ways of partitioning each level of regions but without
stopping in any of the regions. Let $J^{*(k)}$ denote the decision
made during this process in the first $k$ levels of the
recursion. Each realization of $J^{*(k)}$ determines a partition of
$\Omega$ consisting of regions $A\in\mathcal{A}^k$ (not $\mathcal
{A}^{(k)}$ as in the
case of optional stopping). Let $\mathcal{A}^k(J^{*(k)})=\{A\in
\mathcal{A}^k\dvtx A$ is a
region in the partition induced by $J^{*(k)}\}$. If
$A\in\mathcal{A}^k(J^{*(k)})$, then it can be written as
\[
A=\Omega_{l_1l_2\cdots l_k}^{j_1j_2\cdots j_k}.
\]
We set
\[
Q^{*(k)}(A)=\theta_{l_1}^{j_1}\cdot\theta_{l_1l_2}^{j_1j_2}
\cdots\theta_{l_1\cdots l_k}^{j_1\cdots j_k} \quad\mbox{and}\quad
Q^{*(k)}(\cdot|A)=\mu(\cdot|A).
\]
This defines $Q^{*(k)}$ as a random measure.

\item Given the results in Step 1, generate the optional stopping
variables $S=S(A)$ for each region $A\in\mathcal{A}^k(J^{*(k)})$,
successively for each level $k=1,2,3,\ldots.$
Then for each $k$, modify $Q^{*(k)}$ to get $Q^{(k)}$ by replacing
$Q^{*(k)}(\cdot|A)$ with $\mu(\cdot|A)$ for any stopped region $A$ up
to level $k$.
\end{enumerate}
For each $A\in\mathcal{A}^k(J^{*(k)})$, let $I^k(A)=$ indicator of
the event
that $A$ has not been stopped during the first $k$ levels of the
recursion:
\begin{eqnarray*}
E \bigl(Q^{(k)} (T_1^k ) |J^{*(k)} \bigr)
&=& E \biggl(\sum_{A\in\mathcal{A}^k(J^{*(k)})}Q^{*(k)}(A)I^k(A) |J^{*(k)}
\biggr)\\
&=& \sum_{A\in\mathcal{A}^k(J^{*(k)})}E \bigl(Q^{*(k)}(A) |J^{*(k)} \bigr)E \bigl(I^k(A)
|J^{*(k)} \bigr)\\
&\leq& (1-\delta)^k\sum_{A\in\mathcal{A}^k(J^{*(k)})}E \bigl(Q^{*(k)}(A)
|J^{*(k)} \bigr)\\
&=& (1-\delta)^k.
\end{eqnarray*}
Thus $E(Q^{(k)}(T_1^k))\to0$ geometrically and hence
$Q^{(k)}(T_1^k)\to0$ with probability~$1$. Similarly, $\mu(T_1^k)\to0$
with probability $1$.

For any Borel set $B\subset\Omega$, we claim that $\lim Q^{(k)}(B)$
exists with probability $1$. To see this, write
\begin{eqnarray*}
Q^{(k)}(B)&=&Q^{(k)} (B\cap T_0^k )+Q^{(k)} (B\cap T_1^k )\\
&=&a_k+b_k;
\end{eqnarray*}
$a_k$ is increasing since
\begin{eqnarray*}
Q^{(k+1)} (B\cap T_0^{k+1} )&\geq& Q^{(k+1)} (B\cap T_0^k )\\
&=&Q^{(k)} (B\cap T_0^k ),
\end{eqnarray*}
and $b_k\to0$ since $Q^{(k)}(T_1^k)\to0$ with probability $1$.

Since the Borel $\sigma$-field $\mathcal{B}$ is generated by
countably many
rectangles, we have with probability $1$ that $\lim Q^{(k)}(B)$ exists
for all $B\in\mathcal{B}$. Define $Q(B)$ as this limit. If $Q(B)>0$ then
$Q^{(k)}(B)>0$ for some $k$. Since $Q^{(k)}\ll\mu$ by construction, we
must also have $\mu(B)>0$. Thus $Q$ is absolutely continuous.

For any $B\in\mathcal{B}$, $Q^{(k)}(B\cap T_0^k)=Q(B\cap T_0^k)$, and hence
\begin{eqnarray*}
\bigl|Q^{(k)}(B)-Q(B) \bigr|
&=& \bigl|Q^{(k)} (B\cap T_1^k )-Q (B\cap T_1^k ) \bigr|\\
&<& 2Q^{(k)} (T_1^k )\longrightarrow0.
\end{eqnarray*}
Thus the convergence of $Q^{(k)}$ to $Q$ is in variational distance.
\end{pf*}

The next result shows that it is possible to construct optional P\'{o}lya
tree distribution with positive probability on all $L_1$ neighborhoods
of densities.
\begin{theorem}\label{thm2}
Let $\Omega$ be a bounded rectangle in $\mathbb{R}^p$. Suppose that the
condition of Theorem \ref{thm1} holds and that the selection
probabilities $\lambda_i(A)$, the assignment probabilities
$\alpha_i^j(A)/(\sum_l\alpha_l^j(A))$ for all $i,j$ and
$A\in\mathcal{A}^{(\infty)}$ are uniformly bounded away from $0$ and
$1$. Let
$q=dQ/d\mu$; then for any density $f$ and any $\tau>0$, we have
\[
P \biggl(\int|q(x)-f(x)|\,d\mu<\tau\biggr)>0.
\]
\end{theorem}
\begin{pf}
First assume that $f$ is uniformly continuous. Let
\[
\delta(\varepsilon)={\sup_{|x-y|<\varepsilon}}|f(x)-f(y)|;
\]
then $\delta(\varepsilon)\downarrow0$ as $\varepsilon\downarrow0$. For any
$k$ large enough, we can find a partitioning $\Omega=\bigcup_{i=1}^IA_i$
where $A_i\in\mathcal{A}^k$ is arrived at by $k$ steps of recursive
partitioning (deterministic and without stopping) and that each $A_i$
has diameter $<\varepsilon$.

Approximate $f$ by a step function $f^*(x)=\sum_if_i^*I_{A_i}(x),
f_i^*=\int_{A_i}f \,d\mu/\break\mu(A_i)$. Let $D_\varepsilon(f)$ be the set of
step functions $g(\cdot)=\sum g_iI_{A_i}(\cdot)$ satisfying
\[
{\sup_i}|g_i-f_i^*|<\delta(\varepsilon).
\]
Suppose $g\in D_\varepsilon(f)$; then for any $B$ we have
$B=\bigcup_{i=1}^I(B\cap A_i)=\bigcup_{i=1}^IB_i$ and
\begin{eqnarray*}
\biggl|\int_B(g-f)\,d\mu\biggr|&\leq&\sum_i|g_i-f_i^*|\mu(B_i)+\sum_i \biggl|f_i^*\mu
(B_i)-\int_{B_i}f \,d\mu\biggr|\\
&\leq&\sum_i\delta(\varepsilon)\mu(B_i)+\sum_ir_i,
\end{eqnarray*}
where
\begin{eqnarray*}
r_i&=&\mu(B_i) \biggl|\frac{\int_{A_i}f \,d\mu}{\mu(A_i)}-\frac{\int
_{B_i}f \,d\mu
}{\mu(B_i)} \biggr|\\
&=&\mu(B_i) \biggl|\frac{\int_{A_i} (f(x)-f(x_k) )\,d\mu}{\mu(A_i)}-\frac
{\int
_{B_i} (f(x)-f(x_k) )\,d\mu}{\mu(B_i)} \biggr|,
\end{eqnarray*}
where $x_i\in B_i$. Since
\[
|f(x)-f(x_i)|<\delta(\varepsilon) \qquad\mbox{for }x\in A_i,
\]
we have
\[
|r_i|<2\delta(\varepsilon)\mu(B_i).
\]
Hence
\[
\biggl|\int_B(g-f)\,d\mu\biggr|<3\delta(\varepsilon)\mu(B)\qquad \forall B,
\]
and thus
\[
\int|g-f|\,d\mu<3\delta(\varepsilon)\mu(\Omega)=3\delta'(\varepsilon),
\]
where $\delta'(\varepsilon)=\delta(\varepsilon)\mu(\Omega)$. Since all
probabilities in the construction of
$q^k=\frac{dQ^{(k)}}{d\mu}$ are bounded away from $0$ and $1$, we have
\[
P \bigl(q^k\in D_\varepsilon(f)\mbox{ for all large } k \bigr) >0.
\]
Hence
\[
P \biggl(\int|q^k-f|\,d\mu<3\delta'(\varepsilon) \mbox{ for all large } k \biggr)> 0.
\]
On the other hand, by Theorem \ref{thm1}, we have
\[
P \biggl(\int|q^k-q|\,d\mu\to0 \biggr)=1.
\]
Thus
\[
P \biggl(\int|q-f|\,d\mu<4\delta'(\varepsilon) \biggr)>0.
\]
Finally, the result also holds for a discontinuous $f$ since we can
approximate it arbitrarily closely in $L_1$ distance by a uniformly
continuous one.
\end{pf}

It is not difficult to specify $\alpha_i^j(A)$ to satisfy the
assumption of Theorem \ref{thm2}. A~useful choice is
\[
\alpha_i^j(A)=\tau^k\mu(A_i^j ) /\mu(\Omega) \qquad\mbox{for }A\in
\mathcal{A}^k,
\]
where $\tau>0$ is a suitable constant.

The reason for including the factor $\tau^k$ when $A\in\mathcal
{A}^k$ is to
ensure that the strength of information we specified for the
conditional probabilities within $A$ is not diminishing as the depth of
partition $k$ increases. For example, in Example \ref{ex2} each $A$ is
partitioned into two parts of equal volumes; that is,
\[
A=A_1^j\cup A_2^j,\qquad \mu(A_1^j )=\mu(A_2^j )=\tfrac12\mu(A).
\]
Thus $A\in\mathcal{A}^k\Rightarrow\mu(A_i^j)=2^{-(k+1)}\mu(\Omega
)$, and
\[
\alpha_i^j(A)=2^k \frac{\mu(A_i^j)}{\mu(\Omega)}=\frac12 \qquad\mbox{for
all }k.
\]

In this case, by choosing $\tau=2$ we have obtained a nice
``self-similarity'' property for the optional P\'{o}lya tree, in the sense
that the conditional probability measure $Q(\cdot|A)$ will have an
optional P\'{o}lya tree distribution with the same specification for
$\alpha_i^j$'s as in the original optional P\'{o}lya tree distribution for
$Q$.

Furthermore, in this example if we use $\tau=2$ to specify a prior
distribution for Bayesian inference of $Q$, then for any $A\in\mathcal{A}^k$,
the inference for the conditional probability $\theta_1^j(A)$ will
follow a classical binomial Bayesian inference with the Jeffrey's prior
Beta ($\frac12,\frac12$).

\section{Bayesian inference with an optional P\'{o}lya tree prior}\label{sec3}

Suppose we have observed $\mathbf{x}=\{x_1,x_2,\ldots,x_n\}$ where $x_i$'s
are independent draws from a probability measure $Q$, where $Q$ is
assumed to have an optional P\'{o}lya tree as a prior distribution. In
this section we show that the posterior distribution of $Q$ given
$\mathbf{x}$ also follows an optional P\'{o}lya tree distribution.

We denote the prior distribution for $q=\frac{dQ}{d\mu}$ by
$\pi(\cdot)$. For any $A\subset\Omega$, we define
$\mathbf{x}(A)=\{x_i\in\mathbf{x}\dvtx x_i\in A\}$ and $n(A)=\#(\mathbf{x}(A))=$
cardinality of the set $\mathbf{x}(A)$. Let
\[
q(x)=\frac{dQ}{d\mu}(x) \qquad\mbox{for }x\in\Omega
\]
and
\[
q(x|A)=\frac{q(x)}{Q(A)} \qquad\mbox{for }x\in A;
\]
then the likelihood for $\mathbf{x}$ and the marginal density for
$\mathbf{x}$
can be written, respectively, as
\begin{eqnarray*}
P(\mathbf{x}|Q)&=&\prod_{i=1}^nq(x_i)=q(\mathbf{x}),\\
P(\mathbf{x})&=&\int q(\mathbf{x}) \,d\pi(q).
\end{eqnarray*}
The variable $q$ (or $Q$) represents the whole set of random
variables, that is, the stopping variable $S(A)$, the selection variable
$J(A)$ and the condition probability allocation $\theta_i^j(A)$, etc.,
for all regions $A$ generated during the generation of the random
probability measure $Q$.

In what follows, we assume that the stopping rule needed for $Q$ is an
independent stopping rule. By considering how $\Omega$ is partitioned
and how probabilities are assigned to the parts of this partition, we
have
%
%
\begin{equation}\label{eq1}
q(\mathbf{x})=Su(\mathbf{x})+(1-S) \Biggl(\prod_{i=1}^{K^J} (\theta_i^J
)^{n_i^J} \Biggr)q (\mathbf{x}
|\mathbf{N}^J=\mathbf{n}^J ).
\end{equation}
In this expression:
\begin{longlist}
\item$u(\mathbf{x})=\prod_{i=1}^nu(x_i)$ where
$u(x)=\frac1{\mu(\Omega)}$ is the uniform density on $\Omega$.

\item$S=S(\Omega)$ is the stopping variable for $\Omega$.

\item$J$ is the choice of partitioning to use on $\Omega$.

\item$\mathbf{N}^J=(n(\Omega_1^J),\ldots,n(\Omega_{K^J}^J))$
is the
counts of observations in $\mathbf{x}$ falling into each part of the
partition $J$.
\end{longlist}

To understand $q(\mathbf{x}|\mathbf{N}^J=\mathbf{n}^j)$, suppose
$J=j$ specifies a
partition
$\Omega=\Omega_1^j\cup\Omega_2^j\cup\cdots\cup\Omega_{K^j}^j$;
then the
sample $\mathbf{x}$ is partitioned accordingly into subsamples,
\[
\mathbf{x}=\mathbf{x}(\Omega_1^j )\cup\cdots\cup\mathbf{x}(\Omega
_{K^j}^j ).
\]
Under $Q$, if the subsample sizes $n_1^j,\ldots,n_{K^j}^j$ are given,
then the positions of points in $\mathbf{x}(\Omega_i^j)$ within
$\Omega_i^j$
are generated independently of those in the other subregions. Thus
\[
q (\mathbf{x}|\mathbf{N}^J=n^j )=\prod_{i=1}^{K^j}q (\mathbf
{x}(\Omega_i^j ) |\Omega_i^j ),
\]
where
\[
q (\mathbf{x}(\Omega_i^j ) |\Omega_i^j )=\prod_{x\in\mathbf
{x}(\Omega_i^j)}q (x |\Omega_i^j ).
\]
Note that once $J=j$ is given, $q(\cdot|\Omega_i^j)$ is generated
independently as an optional P\'{o}lya tree according to the parameters
$\bolds\rho,\bolds\lambda,\bolds\alpha$ that are relevant within
$\Omega_i^j$. We denote by $\Phi(\Omega_i^j)$ the expectation of
$q(\mathbf{x}(\Omega_i^j)|\Omega_i^j)$ under this induced optional
P\'{o}lya tree
within $\Omega_i^j$.

In fact, for any $A\subset\bigcup_{k=1}^\infty\mathcal{A}^k$, we
have an
induced optional P\'{o}lya tree distribution $\pi_A(q)$ for the
conditional density $q(\cdot|A)$, and we define
\[
\Phi(A)=\int q (\mathbf{x}(A)|A ) \,d\pi_A(q),
\]
if $\mathbf{x}(A)\neq\varnothing$ and $\Phi(A)=1$ if $\mathbf
{x}(A)=\varnothing$.
Similarly, we
define
\[
\Phi_0(A)=u (\mathbf{x}(A)|A )=\prod_{x\in\mathbf{x}(A)}u(x|A)
\]
and $\Phi_0(A)=1$ if $\mathbf{x}(A)=\varnothing$. Note that
$P(\mathbf{x})=\Phi(\Omega)$ and $u(\mathbf{x})=\Phi_0(\Omega)$.

Next, we successively integrate out [w.r.t. $\pi(\cdot)$] the random
variables in the right-hand side of (\ref{eq1}) according to the order
$q(\mathbf{x}|\mathbf{n}^J), \bolds\theta^J, J$ and $S$ (last).
This gives us
%
%
\begin{equation}\label{eq2}
\Phi(\Omega)=\rho\Phi_0(\Omega)+(1-\rho)\sum_{j=1}^M\lambda_j
\frac{D(\mathbf{n}
^j+\bolds\alpha^j)}{D(\bolds\alpha^j)}\prod_{i=1}^{K^j}\Phi
(\Omega_i^j ),
\end{equation}
where $D(\mathbf{t})=\Gamma(t_1)\cdots\Gamma(t_k)/\Gamma(t_1+\cdots+t_k)$.

Similarly, for any $A\in\bigcup_{k=1}^\infty\mathcal{A}^k$ with
$\mathbf{x}(A)\neq\varnothing$, we have
%
%
\begin{equation}\label{eq3}
\Phi(A)=\rho\Phi_0(A)+(1-\rho)\sum_{j=1}^M\lambda_j \frac
{D(\mathbf{n}^j+\bolds\alpha
^j)}{D(\bolds\alpha^j)}\prod_{i=1}^{K^j}\Phi(A_i^j ),
\end{equation}
where $\mathbf{n}^j$ is the vector of counts in the partition
$A=\bigcup_{i=1}^{K^j}A_i^j$, and $M, K^j, \rho, \bolds\lambda^j$,
$\bolds\alpha^j$, etc., all depend on $A$. We note that in the special case
when the choice of splitting variables are nonrandom, a similar
recursion was given in \cite{hutter09}.

We can now read off the posterior distribution of $S=S(\Omega)$ from
equation (\ref{eq2}) by noting that the first term
$\rho\Phi_0(\Omega)$ and the remainder in the right-hand side of
(\ref{eq2}) are,
respectively, the probabilities of the events
\[
\mbox{\{stopped at $\Omega$, generate $\mathbf{x}$ from $u(\cdot
)$\}}
\]
and
\[
\mbox{\{not stopped at $\Omega$, generate $\mathbf{x}$ by one of the $M$
partitions\}.}
\]
Thus $S\sim$ Bernoulli with probability
$\rho\Phi_0(\Omega)/\Phi(\Omega)$. Similarly, the $j$th term in the
sum (over $j$) appearing in the right-hand side of (\ref{eq2}) is the
probability
of the event
\[
\mbox{\{not stopped at $\Omega$, generate $\mathbf{x}$ by using the
$j$th way
to partition $\Omega$\}.}
\]
Hence, conditioning on not stopping at $\Omega$, $J$ takes value $j$
with probability proportional to
\[
\lambda_j \frac{D(\mathbf{n}^j+\bolds\alpha^j)}{D(\bolds\alpha
^j)}\prod_{i=1}^{K^j}\Phi(\Omega
_i^j ).
\]
Finally, given $J=j$, the probabilities assigned to the parts of this
partition are $\bolds\theta^j$ whose posterior distribution is Dirichlet
($\mathbf{n}^j+\bolds\alpha^j$).

By similar reasoning, we can also read off the posterior distribution
of $S=S(A), J=J(A), \bolds\theta^j=\bolds\theta^j(A)$ from (\ref{eq3})
for any $A\subset\mathcal{A}^k$. Thus we have proven the following.
\begin{theorem}\label{thm3}
Suppose $\mathbf{x}=(x_1,\ldots,x_n)$ are independent observations
from $Q$
where $Q$ has a prior distribution $\pi(\cdot)$ that is an optional
P\'{o}lya tree with independent stopping rule, and satisfying the condition
of Theorem \ref{thm2}, the conditional distribution of $Q$ given
$\mathbf{X}=\mathbf{x}$ is also an optional P\'{o}lya tree where, for each
$A\subset
A^\infty$, the parameters are given as follows:
\begin{enumerate}
\item Stopping probability:
\[
\rho(A|\mathbf{x})=\rho(A)\Phi_0(A) /\Phi(A).
\]
\item Selection probabilities:
\[
P(J=j|\mathbf{x})\propto\lambda_j \frac{D(\mathbf{n}^j+\bolds
\alpha^j)}{D(\bolds\alpha^j)}\prod
_{i=1}^{K^j}\Phi(A_i^j ),\qquad j=1,\ldots,M.
\]
\item Allocation of probability to subregions: the probabilities
$\theta_i^j$ for subregion $A_i^j, i=1,\ldots,K^j$ are drawn from
Dirichlet ($\mathbf{n}^j+\bolds\alpha^j$).
\end{enumerate}
In the above, it is understood that $M, K^j, \lambda_j, \mathbf{n}^j,
\bolds\alpha^j$ all depend on $A$.
\end{theorem}

We use the notation $\pi(\cdot|x_1,x_2,\ldots,x_n)$ to
denote this
posterior distribution for~$Q$.

To use Theorem \ref{thm3}, we need to compute $\Phi(A)$ for
$A\in\mathcal{A}^\infty$. This is done by using the recursion (\ref{eq3}),
which says that $\Phi(\cdot)$ is determined for a region $A$ if it is
first determined for all subregions $A_i^j$. By going into subregions
of increasing levels of depth, we will eventually arrive at some
regions having certain simple relations with the sample $\mathbf{x}$.
We can
often derive close form solutions for $\Phi(\cdot)$ for such
``terminal regions'' and hence determine all the parameters in the
specifications of the posterior optional P\'{o}lya tree by a finite
computation. We give two examples.
\begin{example}[($2^p$ contingency table)]\label{ex3}
Let $\Omega=\{1,2\}\times\{1,2\}\times\cdots\times\{1,2\}$ be a
table with
$2^p$ cells. Let $\mathbf{x}=(x_1,x_2,\ldots,x_n)$ be $n$ independent
observations where each $x_i$ falls into one of the $2^p$ cells according
to the cell probabilities $\{q(y)\dvtx y\in\Omega\}$. Assume that $q$
has an
optional P\'{o}lya tree distribution according to the partitioning scheme
in Example \ref{ex1} where $\lambda_j=\frac1{M}$ if there are $M$
variables still available for further splitting of a region $A$, and
$\alpha_i^j=\frac12, i=1,2$. Finally, assume that
$\rho(A)\equiv\rho$ where $\rho\in(0,1)$ is a constant.

In this example, there are three types of terminal regions.
\begin{enumerate}
\item$A$ contains no observation. In this case, $\Phi(A)=1$.
\item$A$ is a single cell (in the $2^p$ table) containing any number
of observations. In this case, $\Phi(A)=1$.
\item$A$ contains exactly one observation, and $A$ is a region where
$M$ of the $p$ variables are still available for splitting. In this
case,
\[
\Phi(A)=r_M=\int q(x) \,d\pi_M(Q),
\]
where $\pi_M(\cdot)$ is the optional P\'{o}lya tree on a $2^M$ table. By
recursion (\ref{eq3}) we have
\begin{eqnarray*}
r_M&=&\rho2^{-M}+(1-\rho)
\Biggl(\frac1{M}\sum_{j=1}^M\frac{B (3/2,1/2 )}{B (
1/2,1/2 )} \Biggr)\cdot r_{M-1}\\
&=&\rho2^{-M}+(1-\rho)\frac12 r_{M-1}\\
&=&\rho2^{-M} \frac{ (1-(1-\rho)^M )}{1-(1-\rho)}+ \biggl(\frac{1-\rho}{2}
\biggr)^M\\
&=&2^{-M}.
\end{eqnarray*}
\end{enumerate}
\end{example}
\begin{example}\label{ex4}
$\Omega$ is a bounded rectangle in $\mathbb{R}^p$ with a partitioning
scheme as
in Example \ref{ex2}. Assume that for each region, one of the $p$
variables is chosen to split it $(\lambda_j\equiv\frac1{p})$, and
that $\alpha_i^j=\frac12, i=1,2$. Assume $\rho(A)$ is a constant,
$\rho\in(0,1)$. In this case, a terminal region $A$ contains either no
observations [then $\Phi(A)=1$] or a single observation $x\in A$. In
the latter case,
\[
\Phi(A)=r_A(x)=\int_Aq(x|A) \,d\pi_A(Q)
\]
and
\begin{eqnarray*}
r_A(x)&=&\frac{\rho}{\mu(A)}+(1-\rho)
\frac1{p}\sum_{j=1}^p\frac{B (3/2,1/2 )}{B (
1/2,1/2 )}\cdot r_{A_{i(x)}^j}(x)\\
&=&\frac{\rho}{\mu(A)}+(1-\rho)\frac12 r_{A_{i(x)}^j}(x),
\end{eqnarray*}
where $i(x)=1$ or $2$ according to whether $x\in A_1^j$ or
$A_2^j$. Since $\mu(A_1^j)=\mu(A_2^j)=\frac12\mu(A)$ for the Lebesgue
measure, we have
\begin{eqnarray*}
r_A(x)&=&\frac{\rho}{\mu(A)}+(1-\rho)\frac12 \biggl[\frac{\rho}{\mu
(A)\cdot
1/2}+(1-\rho)\frac12[\cdots] \biggr]\\
&=&\frac{\rho}{\mu(A)} [1+(1-\rho)+(1-\rho)^2+\cdots]\\
&=&\frac{1}{\mu(A)}.
\end{eqnarray*}
\end{example}
\begin{example}\label{ex5}
$\Omega$ is a bounded rectangle in $\mathbb{R}^p$. At each level, we split
the regions according to just one coordinate variable, according to
a predetermined order; for example, coordinate variable $x_i$ is used to
split all regions at the $k$th step whenever $k\equiv i$ (mod $p$).
In this case, $\Phi(A)$ for terminal regions are determined exactly
as in Example \ref{ex4}. By allowing only one way to split a region,
we sacrifice some flexibility in the resulting partition in
exchange for a great reduction of computational complexity.
\end{example}

Our final result in this section shows that optional P\'{o}lya tree priors
lead to
posterior distributions that are consistent in the weak topology. For
any probability measure $Q_0$ on $\Omega$, a weak neighborhood $U$ of
$Q_0$ is a set of probability measures of the form
\[
U= \biggl\{Q\dvtx\biggl|\int g_i(\cdot) \,dQ-\int g_i(\cdot) \,dQ_0 \biggr|<\varepsilon_i,
i=1,2,\ldots,K \biggr\},
\]
where $g_i(\cdot)$ is a bounded continuous function on $\Omega$.
\begin{theorem}\label{thm4}
Let $x_1,x_2,\ldots$ be independent, identically distributed variables
from a probability measure $Q$, $\pi(\cdot)$ and
$\pi(\cdot|x_1,\ldots,x_n)$ be the prior and posterior distributions
for $Q$ as defined in Theorem \ref{thm3}. Then, for any $Q_0$ with a
bounded density, it holds with $Q_0^{(\infty)}$ probability equal to
$1$ that
\[
\pi(U|x_1,\ldots,x_n)\longrightarrow1
\]
for all weak neighborhoods $U$ of $Q_0$.
\end{theorem}
\begin{pf}
It is a consequence of Schwarz's theorem \cite{schw65} that the
posterior is weakly consistent if the prior has positive probability
in Kullback--Leibler neighborhoods of the true density
\cite{ghosh03}, Theorem 4.4.2. Thus, by the same argument as in
Theorem \ref{thm2}, we only need to show that it is possible to
approximate a bounded density in Kullback--Leibler distance by step
functions on a suitably refined partition.

Let $f$ be a density satisfying $\sup_{x\in\Omega}f(x)\leq M<\infty$.
First assume that $f$ is continuous with modulus of continuity
$\delta(\varepsilon)$. Let $\bigcup_{i=1}^IA_i$ be a recursive
partition of $\Omega$ satisfying $A_i\in\mathcal{A}^k$ and diameter
$(A_i)\leq\varepsilon$. Let
\[
g_i=\sup_{x\in A_i}f(x),\qquad g(x)=\sum_{i=1}^Ig_iI_{A_i}(x)
\]
and $G=\int g(x)\,d\mu$. We claim that as $\varepsilon\to0$, the density
$g/G$ approximates $f$ arbitrarily well in Kullback--Leibler
distance. To see this, note that
\begin{eqnarray*}
0&\leq& G-1=\int(g-f) \,d\mu=\sum_i\int_{A_i} \bigl(g(x)-f(x) \bigr)\, d\mu\\
&\leq&\sum_i\int_{A_i}\delta(\varepsilon) \,d\mu=\delta(\varepsilon)\mu
(\Omega).
\end{eqnarray*}
Hence
\begin{eqnarray*}
0&\leq&\int f\log\bigl(f/(g/G) \bigr) \,d\mu\\
&=&\int f\log(f/g) \,d\mu+\int f\log G \,d\mu\\
&\leq&\log(G)\leq\log\bigl(1+\delta(\varepsilon)\mu(\Omega) \bigr).
\end{eqnarray*}
Finally, if $f$ is not continuous, we can find a set $B\subset\Omega$ with
$\mu(B^c)<\varepsilon'$ such that $f$ is uniformly continuous on $B$. Then
\begin{eqnarray*}
\int(g-f)\, d\mu&=&\int_B(g-f)\, d\mu+\int_{B^c}(g-f)\, d\mu\\
&\leq&\delta(\varepsilon)\mu(\Omega)+M\varepsilon'
\end{eqnarray*}
and the result still holds.
\end{pf}

\section{Density estimation using an optional P\'{o}lya tree prior}
\label{sec:den}
In this section we develop and test the methods
for density estimation using an optional P\'{o}lya tree prior. Two
different strategies are considered. The
first is through computing the posterior mean
density. The
other is a two-stage approach---first learn a fixed tree topology that
is representative of the underlying structure of the distribution, and then
compute a piecewise constant estimate
\textit{conditional} on this tree
topology. Our numerical examples start with the one-dimensional
setting to demonstrate some of the basic properties of optional P\'{o}lya
trees. We then move onto the two-dimensional
setting to provide a flavor of what happens when the dimensionality of
the distribution increases.

\subsection{Computing the mean}
For the purpose of demonstration, we first consider the situation
described in
Example \ref{ex2} with $p=1$ where the state space
is the unit interval and the splitting point of each elementary region
(or tree node) is the middle
point of its range. In this simple scenario, each node has
only one way to divide, so the only decision to
make is whether to stop or not. Each point $x$ in the
state space $\Omega$ belongs to one and only one
elementary region in $A^{k}$ for each $k$. In this case, the posterior mean
density function can be
computed very efficiently using an inductive procedure. (See the
\hyperref[app]{Appendix}
for details.)

In a multi-dimensional setting with multiple ways to split at each
node, the sets in each $A^{k}$ could overlap, and so the computation of the
posterior mean is more difficult. One way to get around this problem
is to place some restriction on how the elementary regions can
split. For example, an alternate splitting rule requires that each
dimension is split in turn (Example \ref{ex5}). This limits
the number of choices to split for each elementary region to one and
effectively reduces the dimensionality of the problem to one.
However, in restricting the ways to divide, one wastes a lot of
computation on
cutting dimensions that need not be cut which affects the variability
of the estimate significantly. We demonstrate this phenomenon in our
later examples.

Another way to compute (or at least approximate) the posterior mean density
is first explored by Hutter \cite{hutter09}. For any point $x \in
\Omega$,
Hutter proposed computing $\Phi(\Omega| x, D)$ and using
$\Phi(\Omega| x, D)/\Phi(\Omega| D)$ as an estimate of
the posterior mean density at $x$. [Here $D$ represents the observed data;
$\Phi(\Omega|D)$ denotes the $\Phi$ computed for the root node given
the observed data points and
$\Phi(\Omega|x,D)$ is computed treating $x$
as an extra data point observed.] This method is general but
computationally
intensive, especially when there are multiple ways to divide each node. Also,
because this method is for
estimating the density at a
specific point, to investigate the entire function one must evaluate
$\Phi(\Omega| x, D)$ on a grid of $x$ values which makes it even more
unattractive computationally. For this reason, in our later
two-dimensional examples we only use the restriction method discussed
above to compute the posterior mean.

\subsection{The hierarchical MAP method}
Another approach for density estimation using an optional P\'{o}lya tree
prior is to
proceed in two steps---first learn a ``good'' partition or tree
topology over the state space, and then estimate the density
conditional on this tree topology. The first step reduces the prior
process from an infinite mixture of
infinite trees to a fixed finite tree. Given such a fixed tree
topology (i.e., whether to stop or not at each step, and if not, which
way to divide), we can easily compute the (conditional)
mean density function. The posterior probability mass over each node is
simply a product of Beta means, and the distribution
within those stopped regions is uniform by construction. So the key
lies in learning a reliable tree
structure. In fact, learning the
tree topology is useful beyond facilitating density
estimation. A representative partition over the state space by itself sheds
light on the underlying structure of the distribution. Such
information is particularly valuable in high-dimensional problems
where direct visualization of the data is difficult.

Because a tree topology depends only on the
decisions to stop and the ways to split, its
posterior probability is determined by the
posterior $\rho$'s and $\lambda$'s. The likelihood of each fixed tree
topology is the product of a sequence of terms in the form, $\rho$,
$1-\rho$, $\lambda_k$, depending on the stopping and splitting
decisions at each node. One seemingly
obvious candidate tree topology for representing the data structure is
the maximum a posteriori (MAP)
topology, that is, the topology with the highest posterior probability.
However, in this setting the MAP topology
often does not produce the most
descriptive partition for the distribution. It biases toward shorter
tree branches in that deeper tree structures simply have more terms
less than 1 to
multiply into their posterior probability. While the data typically
provide strong evidence for the stopping decisions (and so the
posterior $\rho$'s for all but the very deep nodes are either very
close to 1 or very close to 0), this is not the case for the
$\lambda$'s. It occurs often that for an elementary region the data
points are distributed relatively symmetrically in two or more
directions, and
thus the posterior $\lambda$'s for those directions will be much less
than 1. As a consequence, deep tree topologies, even if they
reflect the actual underlying data structure, often have lower
posterior probabilities than shallow trees do. (This failure of
the MAP estimate relates more
generally to the multi-modality of the posterior distribution as well
as the self-similarity of the prior process and deserves more
studies in its own right.)

We propose the construction of the representative
tree topology through a simple top-down sequential procedure. Starting
from the root node,
if the posterior $\rho> 0.5$ then we stop the tree; otherwise we
divide the tree in the direction $k$ that has the highest
$\lambda_k$. (When there is more than one direction with the same
highest $\lambda_k$, the choice among them is arbitrary.) Then we
repeat this procedure for each $A_{k}^{j}$ until all
branches of the tree have been stopped. This can be viewed as a
hierarchical MAP decision procedure---with each MAP decision being made based
on those made in the previous steps. In the context of building
trees, this approach is natural in that it exploits the hierarchy
inherent in the problem.

\subsection{Numerical examples}
Next we apply the optional P\'{o}lya tree prior to several examples of
density estimation in one and two dimensions. We consider the
situation described in
Example \ref{ex2} with $p=1$ and $2$ where the state space
is the unit interval $[0, 1]$ and the unit square $[0,1] \times
[0,1]$, respectively. The cutting point of each coordinate is the middle
point of its range for the corresponding elementary region. For all the optional
P\'{o}lya tree priors used in the
following examples, the prior stopping
probability $\rho=0.5$ and the prior pseudo-count $\alpha=0.5$ for
all elementary regions. The standard P\'{o}lya tree priors examined (as
a comparison) have
quadratically increasing pseudo-counts
$\alpha=\operatorname{depth}^2$ (see \cite{ferg74} and \cite{kraft64}). For numerical
purpose, we stop dividing the nodes if
their support is under a certain threshold which we refer to as the precision
threshold. We used $10^{-6}$ as the precision threshold in the
one-dimensional examples
and $10^{-4}$ in the two-dimensional examples. Note that in the
1D examples, each node has
only one way to divide, and so we can use the inductive
procedure described in the \hyperref[app]{Appendix} to compute the
posterior mean
density function. For the 2D examples, we implemented and tested the
full optional tree as well as a restricted version based on ``alternate
cutting'' (see Example \ref{ex5}).
\begin{example}[(Mixture of two close spiky uniforms)]
\label{den_ex1}
We simulate data from the following mixture of uniforms:
\[
0.5 U(0.23,0.232)+ 0.5 U(0.233, 0.235)
\]
and we apply three
methods to estimate the density function. The first is to compute the
posterior mean density using an optional P\'{o}lya tree prior. The second
is to apply the hierarchical MAP method using an optional
P\'{o}lya tree prior. The third is to compute the posterior mean using a
standard P\'{o}lya
tree prior. The results are
presented in Figure \ref{fig:den_ex1}. Several points can be made
from this figure. (1) A sample size of 500 is
%
%
\begin{figure}

\includegraphics{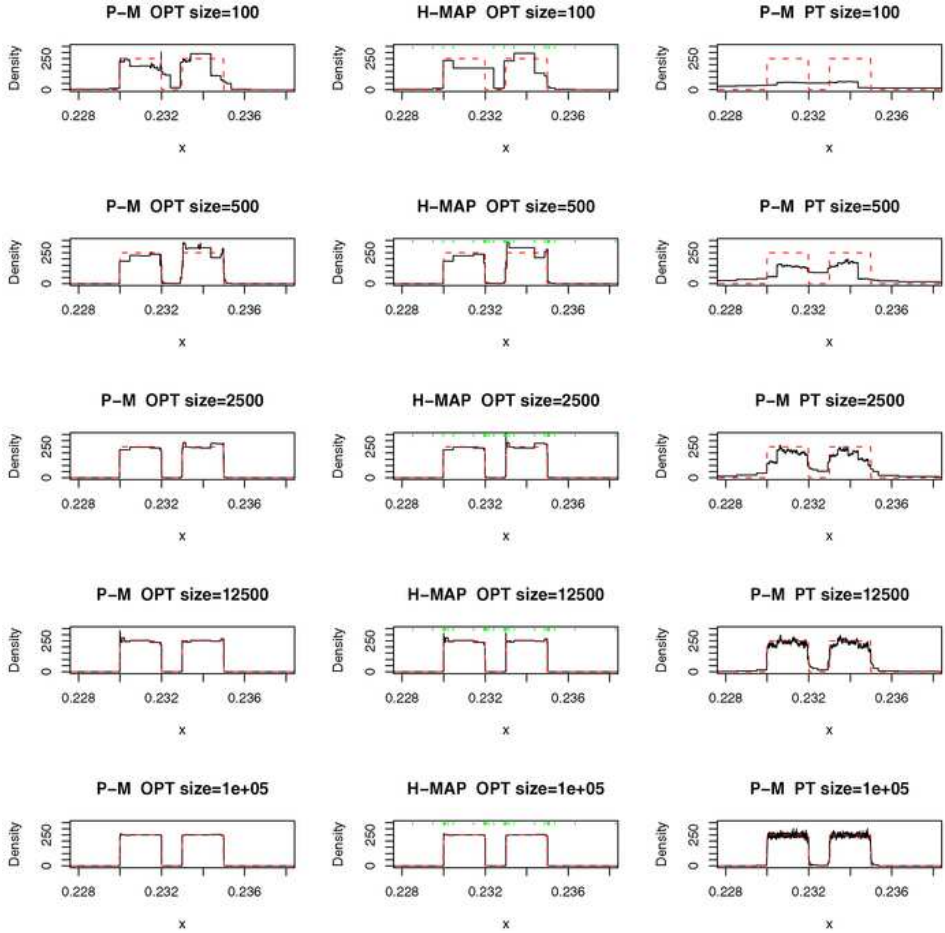}

\caption{Density estimation for $0.5 U(0.23,0.232)+ 0.5
U(0.233, 0.235)$. The five rows represent five different
sample sizes $n= 100$, 500, 2500, 12,500 and 100,000. The first column
corresponds to the posterior mean approach using an optional P\'{o}lya tree
prior. The second column corresponds to
the hierarchical MAP method using an optional P\'{o}lya tree prior. The green
ticks along the top margins of this column
indicate the partition learned from this method. The
third column corresponds to the posterior mean approach using a
standard P\'{o}lya tree prior with $\alpha=\operatorname{depth}^2$. The red dashed
lines in all plots represent the true density function.}
\label{fig:den_ex1}
\end{figure}
sufficient for the optional tree methods to capture the boundaries as
well as the modes of the uniform distributions whereas
the P\'{o}lya tree prior with quadratic pseudo-counts requires thousands of
data points to achieve
this. (2)~With increasing sample size, the estimates from the
optional P\'{o}lya tree methods become smoother, while the estimate from
the standard P\'{o}lya tree
with quadratic pseudo-counts is still ``locally spiky'' even for a
sample size of $10^5$. (This problem can be remedied by increasing the prior
pseudo-counts faster than the quadratic rate at the
price of further loss of flexibility.) (3) The
hierarchical MAP method performs just as well as the posterior
mean approach even though it requires much less computation and memory.
(4) The partition
learned in the hierarchical MAP approach reflects the structure of the
distribution.
\end{example}
\begin{example}[(Mixture of two Betas)]
Next we apply the same three methods to simulated samples from a
mixture of two Beta distributions,
\[
0.7 \operatorname{Beta}(40,60)+0.3
\operatorname{Beta}(2000,1000).
\]
The results are given in Figure
\ref{fig:den_ex2}. Both the
%
%
\begin{figure}

\includegraphics{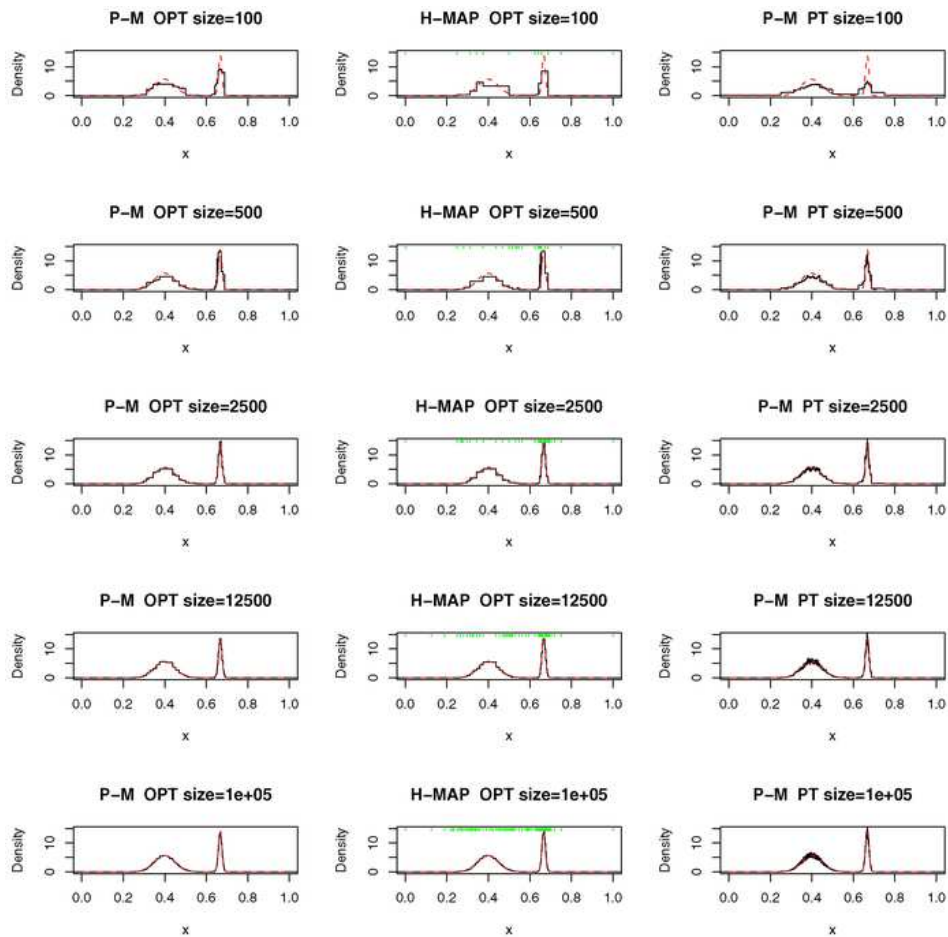}

\caption{Density estimation for $0.7 \operatorname{Beta}(40,60)+0.3 \operatorname{Beta}(2000,1000)$.
The five rows represent five different
sample sizes $n= 100$, 500, 2500, 12,500 and 100,000. The first column
corresponds to the posterior mean approach using an optional P\'{o}lya tree
prior. The second column corresponds to
the hierarchical MAP method using an optional P\'{o}lya tree prior. The green
ticks along the top margins of this column
indicate the partition learned from this method. The
third column corresponds to the posterior mean approach using a
standard P\'{o}lya tree prior with $\alpha=\operatorname{depth}^2$. The red dashed
lines in all plots represent the true density function.}
\label{fig:den_ex2}
\end{figure}
optional and the standard P\'{o}lya tree methods do a decent job in
capturing the locations of the two mixture components (with smooth
boundaries). The optional P\'{o}lya tree does quite well with just 100
data points.
\end{example}
\begin{example}[(Mixture of Uniform and ``semi-Beta'' in the unit square)]
In this example, we consider a mixture distribution over the unit
square $[0,1] \times[0,1]$. The first component is a uniform
distribution over $[0.78, 0.80] \times[0.2, 0.8]$. The second
component has support $[0.25, 0.4] \times[0, 1]$ with $X$ being uniform
over $[0.25, 0.4]$ and $Y$ being
Beta(100, 120), independent of each other. The mixture probability for
the two components is
$(0.35, 0.65)$. Therefore, the actual density function of the
distribution is
\[
\frac{0.35}{0.012} \times{\mathbf{1}}_{[0.78,0.80]\times[0.2,0.8]}
+ \frac
{0.65}{0.15}\times\frac{\Gamma(220)}{\Gamma(120)\Gamma
(100)}y^{99}(1-y)^{119} {\mathbf{1}}_{[0.25,0.4]
\times[0,1]}.
\]
We apply the following methods to estimate this density---(1) the posterior
mean approach using an optional P\'{o}lya tree prior with the alternate
cutting restriction~(Figure \ref{fig:den_ex3_1}); (2) the hierarchical MAP method
using an optional
%
\begin{figure}

\includegraphics{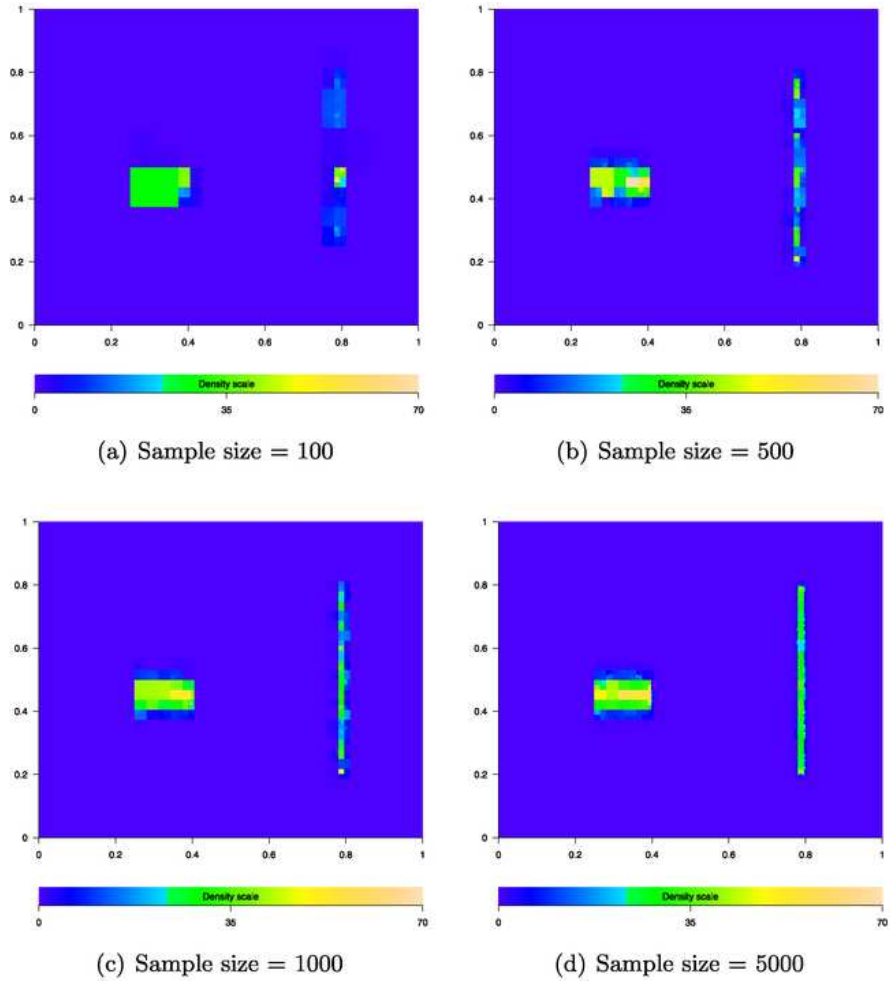}

\caption{Density estimate for a mixture of uniform and
``semi-Beta'' using the posterior mean approach for
an optional P\'{o}lya tree with the restriction of ``alternate
cutting.'' The white blocks represent the density estimates falling
outside of the
intensity range plotted.}
\label{fig:den_ex3_1}
\end{figure}
%
%
\begin{figure}

\includegraphics{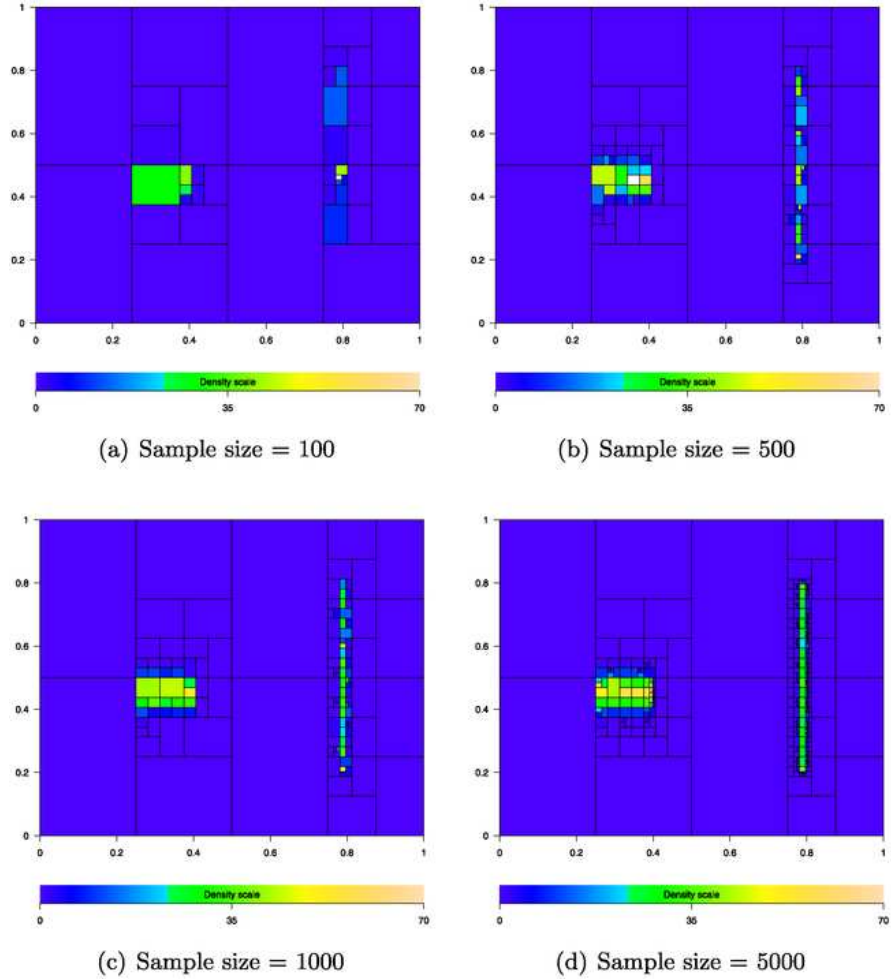}

\caption{Density estimate for a mixture of uniform and ``semi-Beta''
by the hierarchical MAP method using an optional P\'{o}lya tree
prior with the restriction of ``alternate
cutting.'' The dark lines
mark the representative partition learned from the method. The white
blocks represent the density estimates falling outside of the
intensity range plotted.}
\label{fig:den_ex3_3}
\end{figure}
%
%
\begin{figure}

\includegraphics{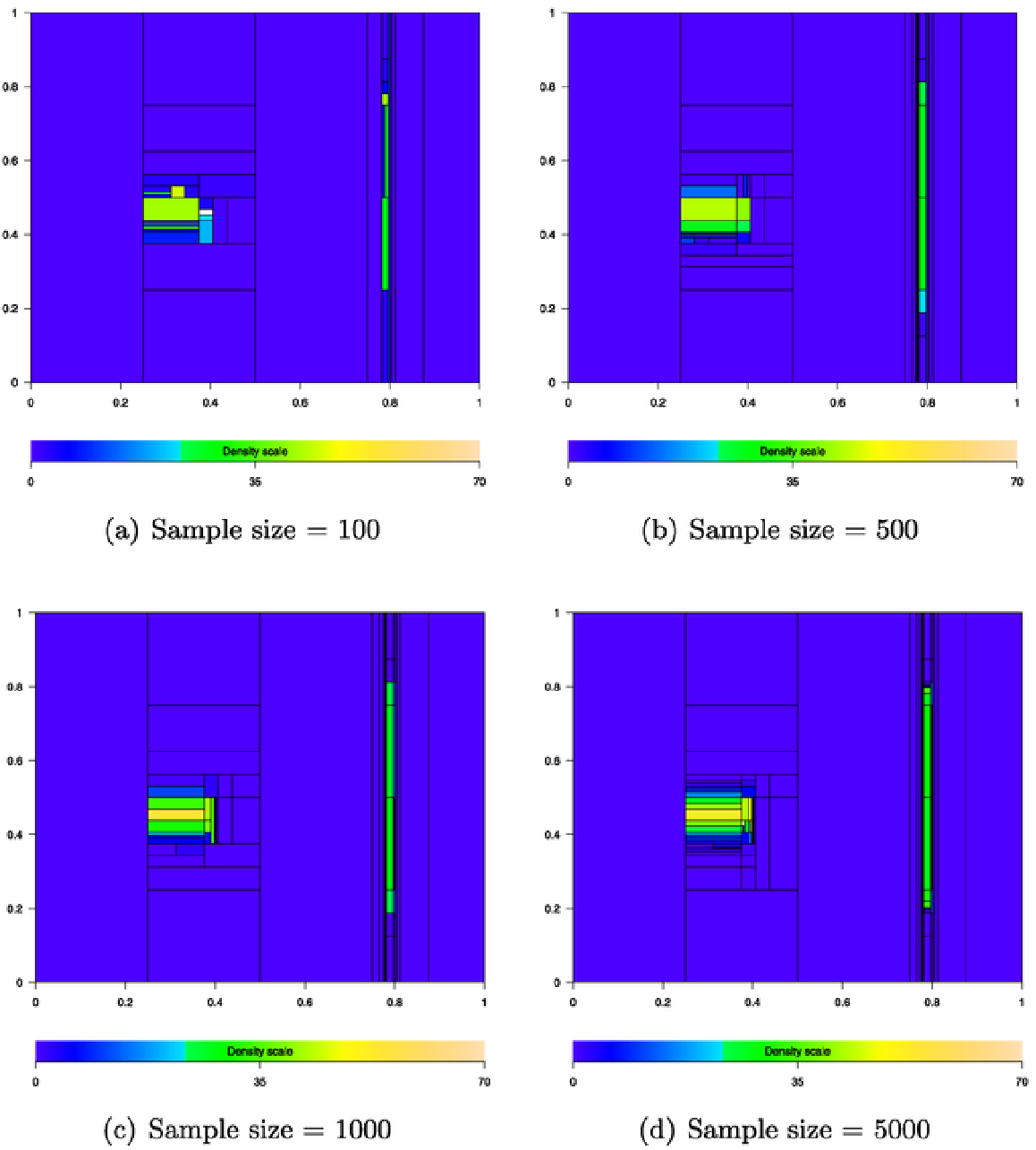}

\caption{Density estimate for a mixture of uniform distribution and
``semi-Beta'' distribution
by the hierarchical MAP method using an optional P\'{o}lya tree
prior (with
no restriction on division). The dark lines
mark the representative partition learned from the method. The white
blocks represent the density estimates falling outside of the
intensity range plotted.}
\label{fig:den_ex3_2}
\end{figure}
P\'{o}lya tree prior with the alternate cutting restriction (Figure \ref{fig:den_ex3_3});
and (3) the hierarchical MAP method using an optional P\'{o}lya tree prior
without any restriction on division (Figure \ref{fig:den_ex3_2}). The last
method does a much better job in capturing the
underlying structure of the data, and thus requires a much smaller sample
size to achieve decent estimates of the density.
\end{example}
\begin{example}[(Bivariate normal)]
In our last example, we apply the hierarchical MAP method using an
optional P\'{o}lya tree prior
to samples from a bivariate normal distribution,
\[
\operatorname{BN}
\left(\pmatrix{0.6 \cr0.4},
\pmatrix{ 0.1^2 & 0 \cr0 &
0.1^2}\right).
\]
This example demonstrates how the posterior optional P\'{o}lya tree
behaves in a multi-dimensional setting
when the underlying distribution has smooth boundary (Figure \ref{fig:den_ex4}).
%
%
\begin{figure}

\includegraphics{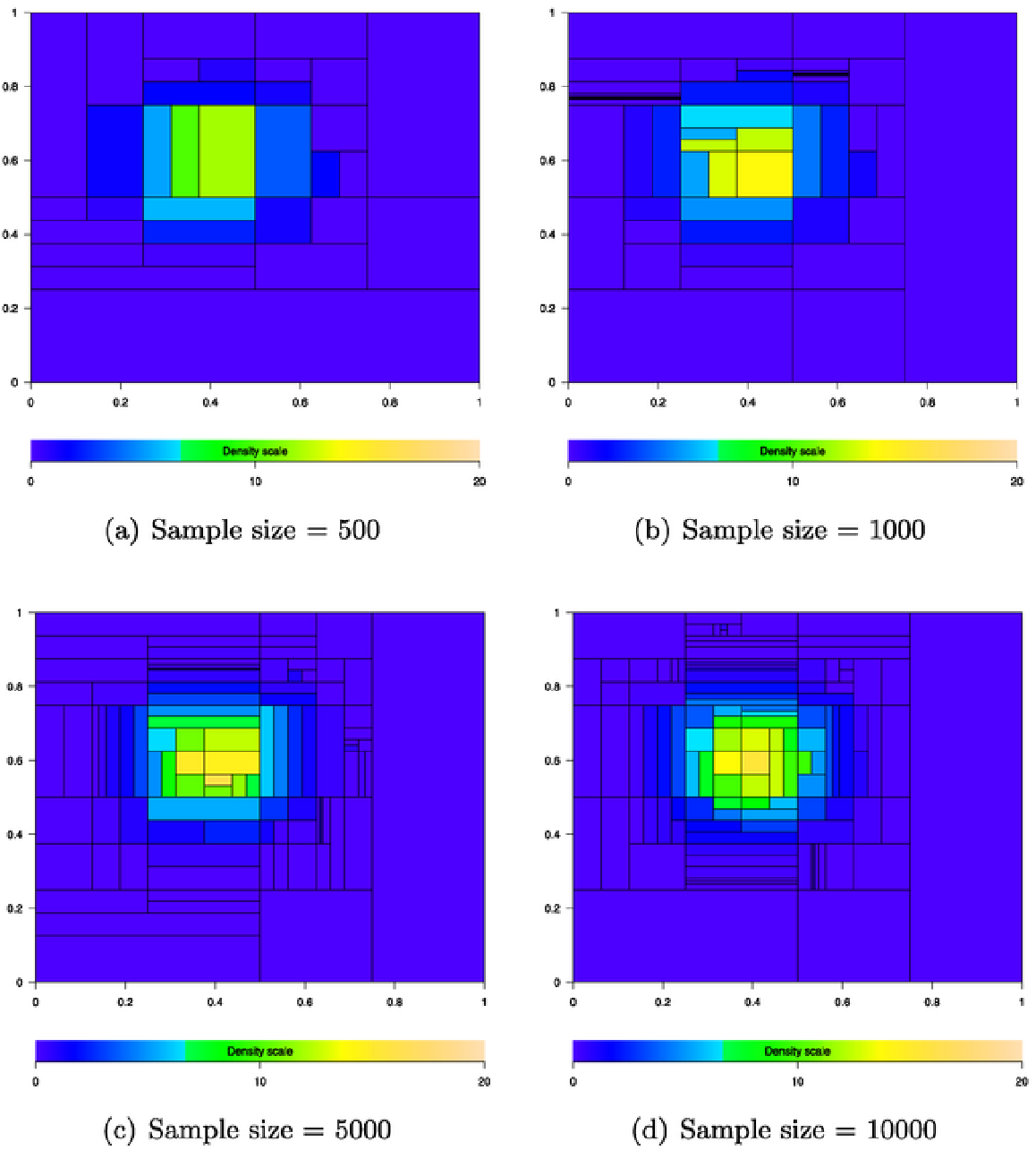}

\caption{The hierarchical MAP method using an optional P\'{o}lya tree prior
applied to samples from a bivariate normal distribution $\operatorname{BN}((0.4,0.6),
0.1^2 I)$.}
\label{fig:den_ex4}
\end{figure}
Not surprisingly, the gradient or change in density is best
captured when its direction is perpendicular to one of the coordinates
(and thus is parallel to the other in the 2D case).
\end{example}

\section{Concluding remarks}
\label{sec:concl}
In this paper we established the existence and the theoretical
properties of absolutely continuous probability measures obtained
through the \hyperref[sec1]{Introduction} of randomized splitting
variables and early
stopping rules into a P\'{o}lya tree construction. For low-dimensional
densities, it is possible to carry out exact computation to obtain
posterior inferences based on this ``optional P\'{o}lya tree'' prior. A
conceptually important feature of this approach is the ability to
learn the partition underlying a piecewise constant density in a
principled manner. Although the theory was motivated by applications
in high-dimensional problems, at present exact computation is too
demanding for such applications. The development of effective
approximate computation should be a priority in future works.

\begin{appendix}\label{app}
\section*{Appendix}

Here we describe an inductive procedure for computing the mean density
function of an optional P\'{o}lya tree when the way to divide each elementary
region is dichotomous and unique.

Let $A_i$ denote a level-$i$ elementary region and $(k_1, k_2,
\ldots, k_i)$ the sequence of
left and right decisions to reach $A_i$ from the root node $\Omega$.
That is,
$A_i = \Omega_{k_1 k_2 \cdots k_i}$, where the $k$'s take values in \{0,
1\} indicating left and right,
respectively.
For simplicity, we let $A_0=\Omega$ represent the root node. Now for
any point $x \in\Omega$, let $\{A_i\}$ be the sequence of
nodes such that $x \in\bigcup_{i=0}^{\infty}A_i$. Assuming $\mu(A_i)
\downarrow0$, the density
of the mean distribution at $x$ is given by
\[
\lim_{i \rightarrow\infty} EP(X\in A_i)/\mu(A_i).
\]
Therefore, to compute the mean density we just need a recipe for
computing $EP(X \in A_i)$ for any elementary region $A_i$. To achieve
this goal, first let $A_{i}'$ be the sibling
of $A_{i}$ for all $i \geq1$. That is,
\[
A_{i}'=\Omega_{k'_1 k'_2 \cdots k'_{i}}\qquad \mbox{where $k'_{j}=k_{j}$
for $j=1, 2, \ldots, i-1$ and $k'_{i}=1-k_{i}$.}
\]
Next, for $i\geq1$, let $\alpha_i$ and $\alpha_i'$ be the Beta
parameters for
node $A_{i-1}$ associated with its two children
$A_i$ and $A_i'$. Also, for $i\geq0$, let $\rho_i$ be the
stopping probability of $A_i$, and $S_i$ the event that
the tree has stopped growing on or before reaching node~$A_{i}$. With
this notation, we have for all $i \geq1$,
\begin{eqnarray*}
&&EP(X\in A_i){\mathbf{1}}(S_i) \\
&&\qquad= EP(X \in A_i){\mathbf{1}}(S_{i-1})
+EP(X \in A_i){\mathbf{1}}
(S_{i-1}^c){\mathbf{1}}(S_{i})\\
&&\qquad= \frac{\mu(A_i)}{\mu(A_{i-1})} EP(X \in A_{i-1}){\mathbf
{1}}(S_{i-1}) \\
&&\qquad\quad{}+
\frac{\alpha_i}{\alpha_i+\alpha_i'}\rho_{i}EP(X \in
A_{i-1}){\mathbf{1}}(S_{i-1}^c)
\end{eqnarray*}
and
\begin{eqnarray*}
EP(X \in A_i){\mathbf{1}}(S_i^c) &=& EP(X \in A^i){\mathbf
{1}}(S_{i}^c){\mathbf{1}}(S_{i-1}^c)\\
&=&\frac{\alpha_i}{\alpha_i+\alpha_i'}(1-\rho_{i})EP(X\in
A_{i-1}){\mathbf{1}}(S_{i-1}^c).
\end{eqnarray*}
Now let $a_i=EP(X\in A_i) {\mathbf{1}}(S_i)$ and $b_i=EP(X \in
A_i) {\mathbf{1}}(S_i^c)$, then the above equations can be rewritten as
%
%
\begin{equation}
\label{eq:ab_den}
\cases{\displaystyle a_i = \frac{\mu(A_i)}{\mu(A_{i-1})} a_{i-1} +
\frac{\alpha_i}{\alpha_i+\alpha_i'} \rho_{i} b_{i-1},\cr
\displaystyle b_i =\frac{\alpha_i}{\alpha_i+\alpha_i'} (1-\rho_{i}) b_{i-1},}
\end{equation}
for all $i\geq1$. Because $a_o=EP(X\in\Omega){\mathbf
{1}}(S_0)=P(S_0)=\rho_0$,
and $b_0=1-a_0=1-\rho_0$, we can apply (\ref{eq:ab_den}) inductively to
compute the
$a_i$ and $b_i$ for all $A_i$'s. Because $EP(X\in A_i)=a_i+b_i$,
the mean density at $x$ is given by
\[
\lim_{i \rightarrow\infty} EP(X\in A_i)/\mu(A_i) = \lim_{i
\rightarrow
\infty} (a_i+b_i)/\mu(A_i).
\]
\end{appendix}\vspace*{-10pt}

\section*{Acknowledgments}

The authors thank Persi Diaconis, Nicholas Johnson and Xiaotong Shen
for helpful comments, and Cindy Kirby for help in typesetting.

\printaddresses


\begin{thebibliography}{99}

\bibitem{black73}
\textsc{Blackwell, D.} (1973).
Discreteness of {F}erguson selections.
\textit{Ann. Statist.} \textbf{1} 356--358.
\MR{0348905}

\bibitem{blackmac73}
\textsc{Blackwell, D.} and \textsc{MacQueen, J. B.} (1973).
Ferguson distributions via {P}\'{o}lya urn schemes.
\textit{Ann. Statist.} \textbf{1} 353--355.
\MR{0362614}

\bibitem{brei84}
\textsc{Breiman, L.}, \textsc{Friedman, J. H.}, \textsc{Olshen, R. A.}
and \textsc{Stone, C. J.} (1984).
\textit{Classification and {R}egression {T}rees}.
Wadsworth Advanced Books and
Software, Belmont, CA.
\MR{0726392}

\bibitem{denison98}
\textsc{Denison, D. G. T.}, \textsc{Mallick, B. K.} and
\textsc{Smith, A. F. M.} (1998).
A {B}ayesian {CART} algorithm.
\textit{Biometrika} \textbf{85} 363--377.
\MR{1649118}

\bibitem{dia86}
\textsc{Diaconis, P.} and \textsc{Freedman, D.} (1986).
On inconsistent {B}ayes estimates of location.
\textit{Ann. Statist.} \textbf{14} 68--87.
\MR{0829556}

\bibitem{fabius64}
\textsc{Fabius, J.} (1964).
Asymptotic behavior of {B}ayes' estimates.
\textit{Ann. Math. Statist.} \textbf{35} 846--856.
\MR{0162325}

\bibitem{ferg73}
\textsc{Ferguson, T. S.} (1973).
A {B}ayesian analysis of some nonparametric problems.
\textit{Ann. Statist.} \textbf{1} 209--230.
\MR{0350949}

\bibitem{ferg74}
\textsc{Ferguson, T. S.} (1974).
Prior distributions on spaces of probability measures.
\textit{Ann. Statist.} \textbf{2} 615--629.
\MR{0438568}

\bibitem{free63}
\textsc{Freedman, D. A.} (1963).
On the asymptotic behavior of {B}ayes' estimates in the discrete
case.
\textit{Ann. Math. Statist.} \textbf{34} 1386--1403.
\MR{0158483}

\bibitem{ghosh03}
\textsc{Ghosh, J. K.} and \textsc{Ramamoorthi, R. V.} (2003).
\textit{{B}ayesian {N}onparametrics}.
Springer, New York.
\MR{1992245}

\bibitem{hans06}
\textsc{Hanson, T. E.} (2006).
Inference for mixtures of finite P\'{o}lya tree models.
\textit{J. Amer. Statist. Assoc.} \textbf{101} 1548--1565.
\MR{2279479}

\bibitem{hutter09}
\textsc{Hutter, M.} (2009).
Exact nonparametric {B}ayesian inference on infinite trees.
Technical Report 0903.5342. Available at \url{http://arxiv.org/abs/0903.5342}.

\bibitem{kraft64}
\textsc{Kraft, C. H.} (1964).
A class of distribution function processes which have derivatives.
\textit{J.~Appl. Probab.} \textbf{1} 385--388.
\MR{0171296}

\bibitem{lav92}
\textsc{Lavine, M.} (1992).
Some aspects of {P}\'{o}lya tree distributions for statistical
modelling.
\textit{Ann. Statist.} \textbf{20} 1222--1235.
\MR{1186248}

\bibitem{lav94}
\textsc{Lavine, M.} (1994).
More aspects of {P}\'{o}lya tree distributions for statistical
modelling.
\textit{Ann. Statist.} \textbf{22} 1161--1176.
\MR{1311970}

\bibitem{lo84}
\textsc{Lo, A. Y.} (1984).
On a class of {B}ayesian nonparametric estimates. {I}. {D}ensity
estimates.
\textit{Ann. Statist.} \textbf{12} 351--357.
\MR{0733519}

\bibitem{maul92}
\textsc{Mauldin, R. D.}, \textsc{Sudderth, W. D.} and
\textsc{Williams, S. C.} (1992).
P\'{o}lya trees and random distributions.
\textit{Ann. Statist.} \textbf{20} 1203--1221.
\MR{1186247}

\bibitem{BarajasandMuller2009}
\textsc{Nieto-Barajas, L. E.} and \textsc{M\"{u}ller, P.} (2009).
Unpublished manuscript.

\bibitem{paddock03}
\textsc{Paddock, S. M.}, \textsc{Ruggeri, F.}, \textsc{Lavine, M.}
and \textsc{West, M.} (2003).
Randomized {P}olya tree models for nonparametric {B}ayesian
inference.
\textit{Statist. Sinica} \textbf{13} 443--460.
\MR{1977736}

\bibitem{schw65}
\textsc{Schwartz, L.} (1965).
On {B}ayes procedures.
\textit{Z. Wahrsch. Verw. Gebiete} \textbf{4}
10--26.
\MR{0184378}

\end{thebibliography}
\end{document}